\documentclass[12pt]{article}
\usepackage{amsmath,amsthm,amsfonts,amssymb,mathrsfs,bm,graphicx,amscd,epsfig,psfrag,verbatim,hyperref}
\usepackage[usenames]{color}
\usepackage{fancyhdr}
\parskip 	\bigskipamount
\newtheorem{theorem}{Theorem}[section]

\newtheorem{remark}{Remark}[section]

\newtheorem{definition}{Definition}

\def\a{{\alpha}}

\def\Z{\mathbb{Z}}

\def\R{\mathbb{R}}
\def\C{\mathbb{C}}
\def\E{\mathbb{E}}
\def\P{\mathbb{P}}

\def\d{\delta}
\def\del{\delta}

\def\la{\lambda}

\def\D{\mathcal{D}}
\def\eps{\epsilon}

\def\b{\mathcal{B}}

\def\ph{\varphi}

\def\var{\mathrm{Var}}

\def\c{\nka}

\def \d{\mathrm{d}}

\def\la{\lambda}
\def\a{\alpha}

\def\tr{\mathrm{tr}}
\def\eps{\varepsilon}

\def\b{V_{\eps}}

\def \del{\delta}

\def \Z{\boldsymbol{Z}}

\def\E{\mathbb{E}}

\def\L{\Lambda}

\def\rt{\rho_\tr^{(2)}}

\def\ph{\varphi}
\def\c{\complement}

\def\eps{\varepsilon}

\def \del{\delta}

\def\b{\beta}

\def\c{\complement}
\def\Var{\mathrm{Var}}
\def\a{\alpha}

\renewcommand{\l}[0]{\left }
\renewcommand{\r}[0]{\right}

\makeatletter
\renewcommand*{\@cite@ofmt}{\hbox}
\makeatother

\begin{document}
%\begin{center}
% \textsc{Notes on rigidity for Coulomb systems}
%\end{center}

\title{Fluctuations, large deviations and rigidity in hyperuniform systems: a brief survey}
\author{
\begin{tabular}{c}
{Subhro Ghosh}\\  
% Princeton University\\ subhrowork@gmail.com
\end{tabular}
\and
\begin{tabular}{c}
{Joel Lebowitz}
% \\ Rutgers University\\ lebowitz@math.rutgers.edu
\end{tabular}
}
\date{}
\maketitle

\begin{center}
\textit{Dedicated to Prof B.V. Rao on the occasion of his 70th birthday}
\end{center}

\begin{abstract}
We present a brief survey of fluctuations and large deviations of particle systems with subextensive growth of the variance. These are called hyperuniform (or superhomogeneous) systems. We then  discuss the relation between hyperuniformity and rigidity. In particular we give sufficient conditions for rigidity of such systems in $d=1,2$.
\end{abstract}

\section{Introduction}
\begin{center}
\textit{ To fluctuate is normal, and in most cases, the fluctuations themselves are normal}
\end{center}
In this brief survey, we explore the subject of fluctuations in several models of hyperuniform particle systems, that is, point processes with reduced number variance. We will also study large deviations for such systems, and finally, the notions of rigidity phenomena in such systems  which has arisen in recent work. 

A quantity of key interest in the study of stochastic particle systems is the fluctuation of the particle number in a domain. More precisely, suppose we have a particle system on a Euclidean space $\R^d$, and suppose we have a sequence of domains $\L_n \uparrow \R^d$ in a self similar manner, that is $\L_n=\{\la_n \cdot x: x \in \L_1\}$ where $0<\la_n \uparrow \infty$. Denoting by $N(\L_n)$ the (random) number of particles in $\L_n$, we are interested in the variance  $\Var(N(\L_n))$. In most models of particle systems, including the Poisson process, Gibbsian models (with tempered interaction potentials), Bosonic and other models exhibiting FKG type properties, the fluctuations are extensive, i.e. asymptotically they grow like the volume: $\Var(N(\L_n))=|\L_n|(1+o(1))$, where $|\L_n|$ denotes the Euclidean volume of $\L_n$. In some cases of physical interest, e.g. at critical points,  they grow faster than $|\L_n|$. When the fluctuations grow like the volume, we call such growth ``extensive''. 

However, there are many  natural models where extensive growth of fluctuations is not true; indeed for thermodynamic limits of Coulomb systems, eigenvalues of random matrices, zeros of random polynomials and many other Fermionic models, the fluctuations are sub-extensive:  $\Var(N(\L_n))=o(|\L_n|)$, and in fact $\Var(N(\L_n))=|\partial \L_n|(1+o(1))$ in many examples. Here $|\partial \L_n|$ denotes the Euclidean area of the boundary $\partial \L_n$ of the domain $\L_n$. Point processes with sub-extensive fluctuations of the particle number are referred to as \textsl{hyperuniform} or \textsl{superhomogeneous}. Hyperuniform processes have been known and studied for several decades (see \cite{MY}, \cite{L}, \cite{Ma}, \cite{ToSt}, \cite{Torq}). Recently they have attracted renewed interest in the material science community (\cite{WiGu},\cite{HeLe}) where hyperuniformity has been claimed in many remarkable contexts like shear flows in dilute suspensions and critical absorbing states in non-equilibrium systems. 
 
Another feature of the particle counts, in a fairly general setting, is that under natural centering and scaling, the fluctuations are asymptotically Gaussian. This is known for a wide range of particle systems (\cite{DV},\cite{L},\cite{Sos2}). Recently, sufficient criteria for the existence of CLT and local CLT, involving the locations of zeros of the generating polynomial for particle count, has been obtained by various authors (\cite{LPRS},\cite{GLiPe}). 

Large deviations (in the space of empirical measures) for particle systems have also been extensively studied (\cite{DZ}, \cite{AGZ}). Other than the case of Gibbsian measures, large deviation results are known for several hyperuniform models, including eigenvalues of Gaussian random matrices and zeros of Gaussian random polynomials. A key instance of this is the study of hole (or overcrowding) probabilities, that is, the event that there are no particles (resp.,  more than typical number of particles) in a large domain. Both moderate and very large deviations are understood (for Gaussian matrices as well as polynomials). These laws are of the same form for both processes (\cite{JLM},\cite{NSV}). 

A relatively recent development has been the study of so-called \textsl{rigidity phenomena}. Roughly speaking, this entails that certain statistics of the particles in a local neighbourhood $\D$ are determined almost surely by the particle configuration outside $\D$. In other words, these statistics of the particles in $\D$ are measurable functions of the particle configuration outside. The most fundamental form of rigidity phenomena is rigidity of the particle number in the domain $\D$.  Following initial results in \cite{GP} and \cite{G}, a wide variety of such rigidity phenomena (and related behaviour) has been studied in a large class of point processes, \cite{Bu}, \cite{BuDQ}, \cite{OsSh}. A very recent result in this direction provides sufficient conditions for rigidity of particle numbers in terms of hyperuniformity and decay of correlations in one and two dimentions \cite{GL}.
 
\section{Basic notions} 

A common general setting in which to study point processes is a locally compact Hausdorff space $X$, equipped with a regular Borel measure $\mu$. We consider the set $\mathcal{S}(X)$ of locally finite point sets on $X$, equipped with the topology of weak convergence on compact sets. It is well known that the space $\mathcal{S}(X)$ is a Polish space with this topology. A point process, formally speaking, is a probability measure on $\mathcal{S}(X)$. Equivalently, it can be seen as a random variable taking values in the space $\mathcal{S}(X)$. Informally, a point process is a random point set in $X$. By identifying a locally finite point set with its induced counting measure, this can also be thought of as a random counting measure on $X$. For a more detailed study of point processes, we refer the reader to \cite{DV}. In this survey, we will mostly specialize to the case $X=\R^d$ and $\mu$ the Lebesgue measure. 

Just as a real-valued random variable is characterized by its cumulative distribution function, similarly the distribution of a point process is described by its various intensity measures. To be precise, the $r$-point intensity measure $\mu_r$ is given by the identity, for $N(D)$ the (random) number of points in any Borel subset $D \subset X$, \[\E\l[ {N(D) \choose r} r! \r] = \int_D \ldots \int_D \d \mu_r(x_1,\ldots,x_r).  \] In most cases $\mu_r$ is absolutely continuous with respect to $\mu^{\otimes r}$, and the corresponding Radon Nikodym derivative $\rho_r$ is called the $r$-point intensity (or correlation) function of the point process. Informally speaking, $\rho_r(x_1,\ldots,x_r)$ denotes the probability density of having points of the process at locations $x_1,\ldots,x_r$. In particular, $\rho_1(x)$ denotes the local particle density per unit measure $\mu$ at $x$, and $\rho_2(x,y)$ denotes the  pair correlation function of the point process. 

For any point process on a Euclidean space $\R^d$, there is a natural way in which a group of translations can act on it. Namely, a translation by a vector $v \in \R^d$ acts on a point configuration $\Upsilon$ as follows: $T_v(\Upsilon):=\{x+v:x \in \Upsilon\}$. Since a point process on $\R^d$ can be thought of as a probability measure on $\mathcal{S}(\R^d)$, therefore this canonically induces an action of the group of translations on a point process. Translation invariance of a point process, therefore, simply means that the law of the point process is invariant under such action. An informal way to understand translation invariance is to say that the statistics of the points in a local neighbourhood does not depend on its location. For a translation invariant point process, all its intensity functions are invariant under the diagonal action of the translation group, and in particular, the one-point intensity function $\rho_1$ is a constant, giving the expected number of particles per unit volume.

In this study, we will  consider point processes on a Euclidean space $\R^d$ that are invariant under the action of the group of translations by $\R^d$ or by $\Z^d$. Unless otherwise stated, our operating assumption will also demand ergodicity of the point process measure under such action. 
%For models on $\Z^d$, we will consider translations only by elements of $\Z^d$, and any domain $\{\L\}$ would be considered as a subset of the lattice $\Z^d$ (and not as a subset of the continuum $\R^d$). 
For periodic models, that is, those models which are invariant in distribution under translations by $\Z^d$, we shall consider the point configuration with a random shift in the unit cube of $\R^d$. This will make the model  invariant under the action of translations of $\R^d$, and will lead to a uniform treatment of the various models under consideration. Key models that we are going to consider will include  the Ginibre ensemble, the Gaussian zero processes, Coulomb systems, determinantal processes and perturbed lattice models. In subsequent sections, we will describe the technical aspects of these models in greater detail.

%The Ginibre ensemble is a special case of Coulomb systems, which are essentially Gibbsian particle systems with a Hamiltonian given by two-body interactions which coincide with the fundamental solution to the Laplacian in the relevant dimension. Coulomb systems in a given dimenson are parametrized by an inverse temperature $\beta$, and the Ginibre ensemble is the two dimensional Coulomb system at inverse temperature $\beta=2$. We will return to a more detailed discussion of Coulomb systems in Section \ref{IntCou}.

\section{Fluctuations in point processes}

\subsection{Fluctuations and hyperuniform processes}
A key object of interest in studying point processes is the particle number. More precisely, for a domain $\L \subset \R^d$, we consider the number $N(\L)$ of points in $\L$. Under our assumptions of translation invariance, it can be easily seen that  in expectation, we have \begin{equation}
\label{mean}
\E[N(\L)]=\rho |\L|,
\end{equation}
where $|\L|$ denotes the Euclidean volume of $\L$, and $\rho$ ($=\rho_1$)  is the (one-point) intensity of the translation invariant point process on $\R^d$.

We can therefore focus our attention on the fluctuations in the particle number. It is known that for ``most'' systems, the size of the fluctuations of $N(\L)$, as measured by their variance $\Var(N(\L))$, will grow like the volume $|\L|$. A typical example is that of a homogeneous Poisson process on $\R^d$. 

Before moving on to the case of sub-volume growth of variance, which will be a key focus in this paper, let us point out that there are examples, particularly in the case of point processes defined on lattices, where we can have $\Var(N(\L))$ grow faster than $|\L|$, i.e. $\Var(N(\L))/|\L| \to \infty$ as $|\L| \uparrow \infty$. Such a phenomenon is observed  at  ``critical points'' in such systems, corresponding to ``critical'' values in the temperature or pressure (\cite{Fi}). 

An important example of such a system is obtained from the Ising spin system with ferromagnetic interactions at zero magnetic field. To map it to a point process, we simply identify the sites having up-spins (or $+$ charges) with having a particle at that site. Under this identification, the variance of $N(\L)$ is $1/4$-th of the variance of the magnetization (which, in turn, is the sum of the signs in the domain $\L$). From classical results on Ising spin systems, it follows that $\var(N(\L))$ grows like the volume $|\L|$ when the inverse temperature $\beta<\beta_c$, where $\beta_c$ is the critical temperature, known to be finite in $d>1$. However, at the critical value of $\beta=\beta_c$, it is known that $\Var(N(\L))$ grows faster than $|\L|$ (in fact, it grows like a power law $|\L|^{\gamma}$ where $\gamma >1 $). For $\beta>\beta_c$, the system is not ergodic, with the variance being extensive in each of the two extremal states. For a detailed reference, we direct the reader to  \cite{Fi}, \cite{Ge}. 

As noted already, our concern here is with hyperuniform systems,  where  the variance is sub-extensive, that is, 
\begin{equation}
\label{subext}
\lim_{\L \uparrow \R^d} \frac{\Var(N(\L))}{|\L|} \to 0.
\end{equation}
%Such processes are referred to in the literature as hyperuniform or superhomogeneous.
%, and have recently attracted considerable interest in probability and statistical physics. These processes are our main objects of interest in this survey.

\subsection{Ginibre's theorem}
Let us begin, however, with an old elegant result by Ginibre (\cite{Gi-1}), providing sufficient conditions for an extensive lower bound on $\Var(N(\L))$, that is, for not being hyperuniform. 
\begin{theorem}[Ginibre]
	\label{Ginibre}
Let $X$ be a random variable taking on integer values in the range $0\le m \le N \le \infty$, with $\P(X=m)=p(m)$. If for some $A>-1$ and all $m \in [0,N-2]$, we have
\begin{equation}\label{gincond} (m+2)\frac{p(m+2)}{p(m+1)} \ge (m+1)\frac{p(m+1)}{p(m)} - A,\end{equation}
then we can conclude that
\[\Var(X) \ge \frac{\E[X]}{1+A}.\]
\end{theorem}

\begin{remark}
	For a translation invariant or periodic point process with $X=N(\L)$ satisfying \eqref{gincond}, this gives $\Var(N(\L))\ge \frac{\rho |\L|}{1+A}$, where $\rho$ is the one-point intensity.
\end{remark}

\begin{proof}
	Here we give a brief sketch of Ginibre's Theorem. To this end, note that \[\sum_{m \ge 0} p(m) [(m+1)\frac{p(m+1)}{p(m)}]=\E[X],\] and 
	\begin{align*}
	& (1+A)^2(\E[X])^2 \\
	= & \l(  \sum  p(m) [(m+1)\frac{p(m+1)}{p(m)} +Am ] \r)^2 \\
	\le & \sum p(m) \l[ (m+1)\frac{p(m+1)}{p(m)} + Am \r]^2.
	\end{align*}
Expanding the squares and using \eqref{gincond} (coupled with the fact that $A>-1$) gives us the conclusion $\Var(X) \ge \frac{\E[X]}{1+A}$, as desired. 
\end{proof}
Ginibre shows (somewhat cryptically) that \eqref{gincond} is satisfied by $X=N(\L)$ for equilibrium systems with tempered potentials (and some hard-core like conditions), thus proving that such systems are not hyperuniform.  This has implications for the nature of phase transitions in such systems, e.g. the density of a fluid is a continuous function of the pressure. More precisely, suppose $\rho$ is the average density and $P$ is the \textsl{pressure} obtained from the grand canonical ensemble in the thermodynamic limit (for details, see \cite{R}, \cite{Ge}).
%More precisely, it can be shown that  the derivative of the density with respect to the pressue is proportional to the ratio of  $\Var(N(\L))$ to the volume $|\L|$.
%In other words, at inverse temperature $\beta$, let $P$ denote the ``pressure'' of the system, that is, \[.\] Let $\rho(=\rho_1)$ denote the one-point intensity (which is a constant due to translation invariance). Needless to say, both $P$ and $\rho$ are functions of $\beta$. 
Then it is known that 
\begin{equation} \label{pressure} \lim_{|\L| \uparrow \infty} \frac{\var(N(\L))}{|\L|} =  \rho \frac{\d \rho}{\d P}.  \end{equation}

If there were to be a discontinuity in the pressure as a function of density (which would correspond to a zeroth order phase transition), then the right hand side in \eqref{pressure} would have to be 0. This would imply that $\var(N(\L))/|\L|$ would have to tend to $ 0$ as $|\L| \to \infty$ : a possibility that is ruled out by Ginibre's theorem. For more details, we refer the reader to \cite{R}, \cite{Ge}, \cite{Fi}.

Ginibre's theorem, in the context of particle systems, explicitly considers Gibbs measures of systems having two body interaction. Ginibre's theorem has been generalized to certain graph counting polynomials that embody many-body interactions; see \cite{LPRS}.

\section{Variance and the pair correlation function}
\label{vapair}
We begin by reminding the reader of some important statistics related to a point process.  For a point process (with intensities absolutely continuous with respect to the Lebesgue measure on $\R^d$), we define the one and two point intensity (or correlation) functions as 
\begin{equation}
\label{onept}
\E[N(\L)]=\int_\L \rho_1(x) \d x
\end{equation}
%We also define the two point intensity function $\rho_2(x,y)$ as 
and
\begin{equation}
\label{twopt}
\E[{N(\L) \choose 2} 2!]=\int \int_{\L \times \L} \rho_2(x,y) \d x \d y
\end{equation}
for all Borel sets $\L \subset \R^d$. 

We also define the truncated pair correlation function $\rt(x,y)$ as 
\begin{equation}
\label{trpair}
\rt(x,y)=\rho_2(x,y)-\rho_1(x)\rho_1(y),
\end{equation}
and the truncated ``full'' pair correlation function $G(x,y)$ as 
\begin{equation}
\label{trfpair}
G(x,y)=\rho_1(x) \del(x,y) + \rho_2(x,y) - \rho_1(x) \rho_1(y),
\end{equation}
where $\del(x,y)$ is the Dirac delta  function.
An equivalent way to understand $G$ is in terms of expectations:
\[ \Var[N(\L)]  = \int_\L \int_\L G(x,y) \d x \d y\]
For translation invariant systems $G(x,y)=G(x-y)$. Observe that for an ergodic translation invariant process, 
$\rho_2(x - y) \to \rho^2$  as $|x-y| \to \infty$, and consequently, $\rt(x-y) \to 0$ and $G(x,y) \to 0$  in that limit.

For a translation invariant system, we have
\begin{align}
\label{trinv}
\begin{split}
&\Var(N(\L)) \\  =&\int \int_{\L \times \L} G(x-y) \d x \d y \\ =& |\L| \int_{\R^d} G(x) \d x - \int_{\R^d} G(x) \a_\L(x) \d x,  
\end{split}
\end{align}
where $\a_\L(x)=\int_{\R^d} \chi_\L((x+y))[1-\chi(y)]\d y$ and $\chi_\L$ is the indicator function of the domain $\L$.

Consider the situation where $|\L| \uparrow \R^d$ in a self-similar way, e.g. by  dilations $\L_R:=\{ R \cdot x : x \in \L  \}$. In such a situation, $\a_\L$ will grow like the surface area $|\partial \L|$ (with $|\partial \L|=2$ for $d=1$). Under mild conditions on $\L$ (e.g. smooth boundaries), $|\partial \L| \sim |\L|^{(d-1)/d}$ as $|\L| \uparrow \R^d$.

Dividing $\Var(N(\L))$ by $|\L|$, we get \[  \lim_{\L \uparrow \R^d} \frac{\Var(N(\L))}{|\L|} = \int_{\R^d} G(x) \d x. \]

\begin{definition}
\textbf{Hyperuniform systems} are those for which \begin{equation} \label{suphom} \lim_{\L \uparrow \R^d} \frac{\Var(N(\L))}{|\L|}= \int_{\R^d} G(x) \d x =0. \end{equation}
\end{definition}
This means that $\int \rt(x)\d x =-\rho$. That in turn implies, in particular, that systems for which $\rt(x)\ge 0$, e.g. those satisfying the FKG inequalities (see \cite{R}, \cite{Ge}), cannot be hyperuniform.

Averaging $\a_\L/|\partial \L|$ over rotations we obtain \begin{equation}  \label{bdry}  \lim_{|\L| \to \infty}  \frac{\a_\L(r)}{|\partial \L|} = \a_d |r| ,\end{equation} where $\a_d$ is a constant (\cite{MY}).

For hyperuniform systems we thus have that the spherically averaged $G(r)$ has the property
\begin{equation}
\label{suphom0}
\int_0^\infty r^{d-1} G(r) \d r = 0
\end{equation}
and
\begin{equation}\label{suphom1} 
\frac{\Var(N(\L))}{|\partial \L|}  = -\a_d \int_0^\infty r^d G(r) \d r \ge 0. 
\end{equation}
In obtaining \eqref{suphom1}, we have combined \eqref{trinv}, \eqref{suphom} and \eqref{bdry}.
$\var(N(\L))$ will grow like $|\partial \L|$ when the right hand side of \eqref{suphom1}, corresponding to the first moment of $G$, exists. 
This implies in particular that $G(r)$ must decay faster than $1/r^{d+1}$. It follows that in $d=1$, bounded variance $\var(N(\L)) \le C$ implies that \[ |\rt(r)| \le \frac{K}{1+r^2}. \] This will be used later in Section \ref{rigidity}.  When the right hand side of \eqref{suphom1} is infinite, $\Var(N(\L))$ will grow faster than $|\partial \L|$ but slower than $|\L|$. This will be the case for the Dyson log gas discussed later.

The question whether $\Var(N(\L))$ can grow slower than $|\partial \L|$ has attracted considerable interest. It was finally settled by Beck in 1987 (see \cite{Be}) where he showed that $\Var(N(\L))$ cannot grow slower than $|\partial \L|$ if the distribution is rotationally invariant (or $\L$ is a ball). It is still an open question as to how slowly this variance can grow, and whether it attains its minimum value for a regular lattice (made translation invariant by averaging over shifts). Interestingly, it has been shown (\cite{CT}, \cite{BC}) that for a simple cubic lattice, there is a transition in some (large enough) dimension $d (\sim 800)$ where putting particles randomly inside each cube gives a smaller variance in a ball than just having particles on $\Z^d$. 

In the translation invariant case, it is relevant to consider the Fourier transform of $G(r)$. Usually denoted as $S(k)$, it  is non-negative, and is referred to as the ``structure function'' in the physics literature (e.g. see \cite{HM}).  This is an important physical quantity in the  study of fluids, where it turns out to be a quantity that can be actually measured experimentally in many situations. It follows from \eqref{suphom} that a  system is hyperuniform when the structure function vanishes at the origin: $S(k) \to 0$ as $|k| \to 0$. A relevant question is how it converges to 0 (as a power law, for example ?)  Such rates are related to the decay of $\rt(r)$ as $r \to \infty$, and thus also, via \eqref{suphom1}, to the  growth of $\var(N(\L))$ in hyperuniform systems. In many physical cases, one expects power law decay: $S(k) \sim |k|^\a$ (as $k \to 0$) and a corresponding decay of $\rt(r) \sim r^{-\gamma}$ (as $r \to \infty$) with $\gamma \ge d+\a$ (where $\a >0$) in order for \eqref{suphom} to hold. For more details, we refer the reader to \cite{ToSt}, \cite{LWL}.
%and we believe that the decay rate of the structure function can be related to various kinds of ``rigidity'' behaviour in particle systems. 

\section{Poisson and other extensive systems}
\label{extensive}
The Poisson point process is the most basic example of a point process; in many ways it is the analogue of the uniform distribution in the world of point processes. A Poisson point process can be defined on any locally compact space $X$ with a background measure $\mu$, and is uniquely characterized by the fact that the points in two disjoint subsets of $X$ are independent of each other, and the one point intensity measure $\mu_1 = \rho \d \mu$. For the homogeneous Poisson point process on $\R^d$ (where \textsl{homogeneous} implies that the background measure $\mu$ is the Lebesgue measure), it is easy to see that the variance of the particle number is \textsl{extensive}. In fact, for the homogeneous Poisson process of intensity $\rho$ and a domain $\L \subset \R^d$, we have the equality  $\Var(N(\L))=\rho|\L|$,  where $|\L|$ denotes the volume of $\L$. 

More generally, we call a particle system ``extensive'' if the following condition is satisfied:  if $\L_n$ is a sequence of domains that are increasing to exhaust $\R^d$ in a self-similar manner, then $\Var(N(\L_n)) \ge |\L_n|(1+o(1))$. Such \textsl{extensive} fluctuations of particle number is also true for for many other systems, including Gibbsian systems with tempered potentials and any non-Gibbsian particle system satisfying the Ginibre Theorem or obeying the FKG inequality (see \cite{FLM}). For the Poisson point process and many systems with extensive variances as well as for some hyperuniform systems, we also have a CLT for the normalized particle number $\l(N(\L_n)- \E[N(\L_n)]\r)/\sqrt{\Var(N(\L_n))}$ (see \cite{CL}, \cite{Ge}).

\section{Coulomb systems}

\subsection{The one component plasma}
Coulomb systems are the primary physical examples of hyperuniform processes.
To simplify matters, we shall consider first the simplest kind of Coulomb system: the classical one component plasma (OCP). This model, also known as ``Jellium'',  was introduced  by Wigner in 1934 \cite{W1}. It consists of particles with a positive charge $e$ moving in a uniform background of negative charge with density $-\rho e$. The background produces an external potential  proportional to $\rho e r^{2}_{i}$; where $r_i$ is the distance of the $i$-th particle from the center of rotational symmetry. This model, as we shall see later, is also of interest in other contexts, such as the distribution of eigenvalues of random matrices.

Setting $e=1$, the potential energy of such a system of $N$ particles in a spherical domain in $\R^d$ (or the whole of $\R^d$) is given by 
\begin{equation} 
\label{pot en}
U(x_1,\cdots,x_N)= \sum_{i<j}^N v_d(x_i-x_j) + \frac{\rho}{2} \sum_{i=1}^N |x_i|^2,
\end{equation}
where, setting $r=|x_i-x_j|$, we have \[  v_d(r)= \begin{cases}   -r &\mbox{if } d=1 \\ -\log r &\mbox{if } d=2   \\  r^{2-d} &\mbox{if } d \ge 3. \end{cases} \]
One can also consider this system in a periodic box or on the surface of a sphere (by setting $v_d(x)=\sum_{m=-\infty, \ne 0}^{\infty} \frac{1}{m^2} \exp[-2 \pi m x / L]$), see  \cite{Ma}. 
 
 The canonical equilibrium probability distribution of this system is given by
\begin{equation} 
\label{eqdist}
\mu_N \propto \exp[-\beta U].
\end{equation}
When $N \to \infty$, the measures $\mu_N$  are expected (and proven in some cases)  to have  a limit $\mu$, which describes a random point process in $\R^d$ with average particle density $\rho$. The extremal measures of the limiting process are (expected to be) translation invariant or periodic (\cite{BM}, \cite{Im}).

This system is exactly solvable in $d=1$: the extremal $\mu$ is periodic with period $\rho^{-1}$, for all $\beta>0$, (see \cite{Ku} and \cite{AM}). The probability distribution of $(N(\L) - \rho |\L|)$, $\L$ an interval, has exponential decay with an exponent that has a nonzero limit as $|\L| \to \infty$ (\cite{MY}). The variance is therefore bounded, and is trivially proportional to $|\partial \L|=2$. This is an example of the general fact that  extremal measures for general 1D systems with bounded variance (or at least tightness of $N(\L) - \rho |\L|$) are periodic (\cite{AGL}).

In $d \ge 2$, the system is translation invariant at  ``small'' $\beta$. For ``large'' $\beta$, the system is expected to form a periodic ``Wigner crystal''.  Numerical simulations predict the formation of the Wigner crystal to be around $\beta =140$, in $d=2$. In $d=2$ this system is exactly solvable at $\beta =2$, where it has the same distribution as the eigenvalues of an i.i.d. complex Gaussian matrix, namely the Ginibre ensemble, scaled to have average density $\rho$. The Ginibre ensemble was introduced by  J. Ginibre as a non-Hermitian analogue of Wigner's Hermitian random matrix models for complex Hamiltonians \cite{Gi-2}. In particular, one has an exact expression for the correlation functions, which have excellent clustering properties, with the truncated pair correlation functions decaying like a Gaussian (\cite{J2}) : 
\begin{equation} \rho_2(r) - \rho^2 = -\rho^2 e^{-\pi \rho r^2}, r=|x_1-x_2|. \label{Tpair}\end{equation} Higher order truncated correlations also decay like $e^{-\gamma D^2}$, where $D$ is the distance between groups of particles. Integrating Eq.~(\ref{Tpair}), one sees that $\int_0^\infty G(r) dr =0$, so this system is hyperuniform.  This is expected to be true for all values of $\beta$ and all $d$ due to Debye screening of charges (\cite{Ma}). 

\subsection{Multi-component Coulomb systems}

In multi-component Coulomb systems, we have natural extensions of the various correlation functions. More specifically, suppose there are two species of particles, denoted by  $\a$ and $\b$. Instead of one and two particle intensities $\rho_1$ and $\rho_2$, we have two types of one-particle densities $\rho_\a$ and $\rho_\b$, and three types of  two-particle densities, denoted $\rho_{\a,\a},\rho_{\b,\b},\rho_{\a,\b}$. If $e_\a$ is the charge corresponding to the particles of type $\a$, then we can consider the one-particle charge  intensity, $q_1(x)=\sum_\gamma e_\gamma \rho_\gamma(x)$ and the charged truncated two-particle density 
\[ q_2^{\tr}(x,y)=\sum_{\gamma,\la} e_{\gamma} e_{\la} [\rho_{\gamma,\la}(x,y)-\rho_\gamma(x)\rho_\la(y)]. \] As an analogue of $N(\L)$, we consider the net charge (i.e., the sum total of the charges of the different kinds of particles) $Q(\L)$ in a domain $\L$.
We then have \[ \E[Q(\L)]= \int q_1(x) \d x \] and \[ \Var[Q(\L)] = \int_\L \int_\L [q_1(x) \del (x-y) +  q_2^{\tr}(x,y)] \d x \d y.  \]
For neutral translation invariant Coulomb systems, we have \[ q_1(x) \equiv 0,  \]  and \begin{equation} \label{multiref}  \lim_{\L \uparrow \R^d} \frac{\var(Q(\L))}{|\L|} = \int q_2^{\tr}(x) \d x = 0.  \end{equation}
%We can proceed to define the various quantities like the (charged) truncated total correlation function $G_Q$, and t

The fluctuations in multi-component Coulomb systems are those of the net charge $Q_\L$ (see \cite{Ma} and the references therein). This is in analogy to the fluctuations of $N(\L)$ in the OCP. 
The arguments in Section \ref{vapair} regarding hyperuniformity would go through in this more general setting. The consequences thereof, including rigidity  also follow from similar arguments. 

We note that one may consider variances of any combination of particle numbers of different species in any multi-component system. The definitions of $q_1$ and $q_2$ would be as above, with the  $e_\gamma$ being arbitrary real weights instead of physical charges.

The basic physical reason for this reduction in charge fluctuations in Coulomb systems is the long range nature of the Coulomb force. This causes \textit{shielding} of bare charges by ``Debye screening''.  This means that if there is a fixed charge at the origin, the other charges will arrange themselves in such a way that the electric field produced by the charge is canceled. Mathematically, it was shown by many authors in the 70's and 80's that shielding is a necessary condition for having at least some kind of of clustering of correlation functions (\cite{Ma}). This screening leads to a whole series of ``sum-rules'', of which \eqref{multiref}  is the first one. For details we direct the reader to the  reviews \cite{Ma} and \cite{BM}.

We note that in many physical situations, such as those involving fluids at low and moderate termperatures, we usually consider macroscopic systems as made up of neutral atoms or molecules interacting via effective short range potentials. In such cases, the flcutuations in the net charge $Q(\L)$ in a region $\L$ will be due entirely to the surface of $\L$ cutting these entities in a ``random'' way. $\Var[Q_\L]$ may then be expected to be proportional to the surface area $|\partial \L|$ of $\L$ (\cite{MY},\cite{Ma}).

\subsection{Asymptotic Normality}
For charge-neutral and translation invariant Coulomb systems in $d \geq 2$ the charge fluctuations satisfy a central limit theorem : deviation from the average divided by the square root of the variance gives
		\begin{equation*}
			\frac{Q(\Lambda) }{\sqrt{\Var(Q(\Lambda))}} \rightarrow \xi,
		\end{equation*}
	a standard Gaussian random variable (\cite{MY}).  In fact, if $\Var[Q(\L)] \sim |\partial \L|$, a joint central limit type behaviour is true in the following sense (\cite{L}).
Let $\R^d, d \ge 2$ be divided into cubes $\Gamma_{j}$ of volume $L^{d}$ whose centers are located at the sites $L\mathbb{Z}^{d}$.
Set
$$ \Upsilon_{j} = Q(\Gamma_j)/\sqrt{\Var(Q(\Gamma_j))}$$
The joint distribution of the $\{ \Upsilon_{j} \}$ approaches as $L \to \infty$ a Gaussian measure with covariance
\begin{align*}
  C_{j,k} =  \left[\delta_{j,k} - \frac{1}{2d} \sum_{e}\delta_{j-k, e}\right] = \frac{1}{2d} \left[-\Delta\right]_{j,k}, \tag{$*$}
\end{align*}
where $e$ is the unit lattice vector and $\Delta$ is the discrete Laplacian. This means that the charge fluctuations in $\Gamma_{j}$ are compensated by the opposite charges in neighboring cubes.  
This is exactly what one would expect when the charges are bound together in neutral molecules. 

\section{Determinantal processes}
Determinantal processes are ones for which the $k$-point correlation $\rho_k(x_1,\dots,x_k)=\det[K(x_j,x_l)]_{j,l=1,\dots,k}$. $K$ is Hermitian and all its eigenvalues $\lambda_j$ are in $[0,1]$. There are more general determinantal processes but we shall not consider them here. Determinantal point processes whose kernels are projection operators are hyperuniform (c.f. Soshnikov, \cite{Sos1}).

Key examples of determinantal processes include distribution of eigenvalues of the Ginibre ensemble, which, as already stated, is the same as the 2D OCP at inverse temperature $\beta=2$. It also includes 1D bulk eigenvalue limit of the Gaussian or the Circular Unitary ensembles, a.k.a. the sine kernel process or the Dyson log gas. This also turns out to be a Coulomb system with 2D logarithmic interactions, confined to a line, at inverse temperature $\beta=2$. In this case, $G(r)$ decays like $r^{-2}$ so its first moment is infinite and the variance of the particle number in an interval of length $|\L|$ grows like $\log |\L|$. The Dyson log gas is hyperuniform for all $\beta$ (\cite{For}).   The ground state of an ideal Fermi gas in any dimension is also known to be a determinantal process with a projection kernel, and thus hyperuniform.

One can prove for all determinantal processes a local CLT, using the fact that the zeros of the generating function of a determinantal point process (whether projection or not) all lie on the negative real axis on the complex plane, \cite{ForL} and \cite{CL}.

\section{Perturbed lattice models}
\label{lattice}
 We consider I.I.D. perturbations of a lattice, i.e. each lattice point $z \in \Z^d$ is shifted to $z+x\in \R^d$ with a probability distribution $h(x)dx$. These are like displacements of atoms in an ideal crystal. The resulting processes are (periodic) hyperuniform. This can be seen by noting that the (periodic) one particle density is given by \[ \rho_1(x)=\sum_{z\in \Z^d}h(x-z), \int h(x) dx=1,\] and \[G(x,y) = \rho_1(x)\delta(x-y)-\sum_{z\in\Z^d}h(x-z)h(y-z), \textrm{ so} \int G(x,y)dy=0.\] These systems have $\Var(N_{\Lambda})\sim c|\partial \Lambda|$ when the first moment of $h$ exists and thus bounded variance in 1D (\cite{GaSz}).

\section{G processes}
 Various examples of perturbed lattice models in 1D with bounded variance have been studied in the statistics literature.
 A related model, the G process, was studied in \cite{GLS} as a statistical mechanical point process.
 To construct this process, we define a real-valued Markov process $Y_{\lambda}(t)$, for $t\ge0$, satisfying $Y_{\lambda}(t)>-1$; here $\lambda$ is a probability measure on $(-1,\infty)$.
  
 $Y_{\lambda}(t)$ is defined by two conditions:

 (1) $Y_{\lambda}(0)$ is distributed according to $\lambda$, and

 (2) $Y_{\lambda}(t)$ increases at rate 1 as $t$ increases, except at points of a Poisson process of density $\alpha$ on $\mathbb{R}_+$, at which it jumps down by one unit -- unless this jump would violate the condition $Y_{\lambda}>-1$, in which case no jump occurs.
 
This process has a unique stationary single-time distribution $\lambda=\lambda_0$.
The corresponding translation invariant process (obtained e.g. by imposing the initial condition $\lambda_0$ at time $\tau$ and then taking the Cesaro limit as $\tau \to -\infty$) is denoted by $Y(t)$. The points of the G process are those points at which $Y$ jumps. In other words, the G process is the distribution of the jump points of the $Y$ process. The points of the G process may be viewed as the output of a so-called D/M/1 queue. It is shown in \cite{GLS} that for this process with $\a>1, \rho=1$, $\mathrm{Var}(N_{\L})\le \mathrm{const.}$, for $\L=[s,t]$. It is also shown that this system has exponential decay of the (truncated) pair correlation function.

\section{Gaussian Zeros}
Another important class of hyperuniform point processes that we will consider on $\R^2$ are the zeros of the so-called planar Gaussian analytic function. These are large $N$ limits of the zeros of random polynomials. The \textsl{standard} planar Gaussian zero process is the large $N$ limit of the zeros of the \textsl{Weyl polynomials}, given by \[ \mathfrak{p}_N(z)=\sum_{k=0}^N \xi_k \frac{z^k}{\sqrt{k!}}. \] This is a special case of  the $\a$-Gaussian zeros, which are large $N$ limits of the zeros of $\a$-Weyl polynomials \[ \mathfrak{p}_N^{(\a)}(z)=\sum_{k=0}^N \xi_k \frac{z^k}{(k!)^{\a/2}}. \] 
%Notice that the standard Weyl polynomials are a special case with $\a=1$. 
Like the Ginibre eigenvalues and the Coulomb systems, the standard Weyl polynomials also originate in physics, and have been studied extensively by Bogomolny, Bohigas, Lebeouf and others in the context of spectral analysis of Hamiltonians of chaotic quantum systems (\cite{BBL1},\cite{BBL2}). The $\a$-Gaussian zeros are known to be hyperuniform for $\a>0$.

In extensive work by Nazarov, Tsirelson, Sodin and others (\cite{ST}, \cite{NS1}, \cite{NS2}, \cite{NSV}), it has been shown that the standard planar Gaussian zero process, like the Ginibre ensemble, exhibits translation invariance and Gaussian decay of the truncated pair correlation function. The fluctuations of the particle number are sub-extensive : in fact, we have $\Var(N(\L)) \sim |\partial \L|$ as $\L \uparrow \R^2$ in a self similar manner. Such similarities in behaviour with the Ginibre ensemble calls for a comparative study of the Gaussian zeros and the Ginibre ensemble, and we will see that in spite of the striking similarities between the two, there are spectacular differences as stochastic processes, particularly in the light of \textsl{rigidity phenomena}.

\section{Large deviations}
 As might be expected from the reduction of fluctuations, the probability of large deviations
 from the mean will be smaller for hyperuniform systems than those for systems with Poisson-type fluctuations. This problem was studied for Coulomb systems in \cite{JLM}, using electrostatic type arguments. They found that this is indeed the case in all dimensions and all $\beta>0$.

For the 2D OCP with density $\rho$, the probability of 
having $n(R)$ particles in a disc of radius $R$, corresponding
to a charge   $|Q|= |n(R)-\pi \rho R^2|$, behaves as
$$ \text{Prob} \left\{ | n(R) - \rho \pi R^{2} | > b_{\alpha} R^{\alpha} \right\} \sim \exp\left[-c_{\alpha} R^{\phi(\alpha)}\right], $$
with
$$ \phi(\alpha) = \left\{ \begin{array}{l c r}
  2\alpha - 1 &, & \frac{1}{2} < \alpha \le 1 \\
	3\alpha -2 &,& 1 \le \alpha \le 2 \\
	2\alpha &,& \alpha \ge 2.
\end{array} \right. $$
 
The situation in $d=3$ is similar to that in $d=2$ although the details differ.
 
This probability is much smaller than the large deviations for systems with
short range interactions where, e.g. for $\alpha = 2$ one would get $e^{-cR^{2}}$ instead of $e^{-cR^{4}}$. The symbol $\sim$ means that taking the logarithm of both sides and dividing by $R^{\phi(\alpha)}$ we get a finite limit when $R \rightarrow \infty$. 

These ``macroscopic'' results can be checked and confirmed at $\beta=2$ where we have explicit
solutions for the correlation functions. We can get then additional information such as the
charge density outside a disc of radius $R$ conditioned on there being no particles inside.
In particular the density at $r=R^{+}$ is given by $\rho(R^{+}) \sim \frac{1}{2} \pi \rho^2 R$.

It turns out that the large deviation function for the 2D OCP is of the same 
form, in its dependence on $\alpha$ as that of the point process generated by the zeroes of the standard planar Gaussian Analytic Function (henceforth GAF), $f=\sum_{k=0}^\infty \frac{\xi_{k}}{\sqrt{k!}} z^{k}$, with the $\xi_{k}$ i.i.d standard complex Gaussians (\cite{NSV}). 

 For $d=1$ with $v_1(r)$ (linear) Coulombic interactions, we have already noted the the variance of particle numbers remains bounded in the size of the interval. The probability \[\Pr\{|N(L)-\rho L|>K\}\sim\exp[-cK],\] in any interval of length $L$. Large deviations for this system are expected to behave as (\cite{For})   \[\Pr\{|N(L)-\rho L|>\kappa L\}\sim \exp[-cL^3].\] For $d=1$, with $v_2(r)=-\log r$ interactions \[\Pr\{|N(L)-\rho L|>bL\}\sim\exp[-cL^2].\] For perturbed lattice systems \[ \Pr\{|N(L)-\rho L|>bL\}\sim h(L)^{cL}.\]  

On the other hand, for $G$ processes, this probability goes like $\exp[-cL]$ (see \cite{GLS}).

%\section{Rigidity phenoemena}
\section{Spatial conditioning and DLR equations}
So far we have discussed fluctuations and large deviations of particles, or charges, in a region $\Lambda$ without saying anything about the configuration of particles/charges outside $\Lambda$, i.e. in $\Lambda^{c}=\R^d \setminus \Lambda$. We ask now: what can we say about the distribution of points inside $\Lambda$ given the configuration in $\Lambda^c$, i.e, we want the conditional probability  $\mu\left( dX_{\Lambda}|X_{\Lambda^{c}} \right)$  of a configuration in $dX_{\Lambda}$ given $X_{\Lambda^{c}}$.

For equilibrium Gibbs measures $\mu$  of particle systems on $\R^d$ the answer to this is given by the Dobrushin-Lanford-Ruelle (DLR) equations \cite{R}). 

  \begin{equation}
    \mu \left( x_1, \ldots, x_N | X_{\Lambda^c} \right) = \frac{\mbox{ exp }[-\beta U (X_{\Lambda}|X_{\Lambda^c})] }{ \int e^{-\beta U (X_{\Lambda}| X_{\Lambda^{c}}  ) } d X_{\Lambda} }
    \label{eq:DLR}
  \end{equation}
where  $U(X_\Lambda|X_{\Lambda^c})$ is the potential energy of a configuration in $\Lambda$ given the configuration in $\Lambda^c = \R^d\setminus \Lambda$.

When the interaction $U$ decays sufficiently rapidly with distance and $\mu$ is ergodic, the behaviour of $\Var[N(\Lambda)]$, for large $\Lambda$, is similar to the unconditional case, and the Ginibre lower bound on the variance holds.  This is however not the case for systems with long range Coulomb interactions, where $U(X_\L|X_{\L^c})$ is not well defined. In that case, as we have seen before, the condition for the Ginibre Theorem does not hold, and there  is no strictly positive lower bound on $\var[N(\L)|X_{\L^\c}]$.

\section{Number Rigidity}
\label{rigidity}
The property that the measure $\P(N(\L)|X_{\L^c})$ is concentrated at a single value of $N(\L)$ has been called ``[number] rigidity'' in \cite{GP}. They showed that the Ginibre ensemble and the standard planar Gaussian zero process have this property.  In \cite{G}  number rigidity was also shown for the GUE (and the CUE) point processes. Both the Ginibre and the GUE ensemble correspond to, as already mentioned, Coulomb systems (with logarithmic interactions) at particular temperatures.

\cite{GP} also showed that while $N_{\Lambda}$ is fixed by $X_{\Lambda^c}$, the distribution of points inside $\Lambda$  is not rigid; in fact it is absolutely continuous with respect to the Lebesgue measure. A similar behaviour is true for the $d=1$ Coulomb system considered in \cite{AM}. There it was proved, for $d=1$ Coulomb systems, that the charge in an interval $[a,b]=\L$,  which 
corresponds for the OCP to the number of particles in $\Lambda$, is 
uniquely specified by the configuration $X_{\Lambda^{c}}$ for \textit{all typical} configurations with respect to the infinite volume measure $\mu$. (The set of atypical configurations has measure zero). 

After the work of \cite{GP} and \cite{G}, various authors have established rigidity for a number of point processes, e.g. Beta, Gamma and Airy processes (\cite{Bu}). In all these cases, the process for which rigidity was proven is hyperuniform. 

In \cite{GL} it has been shown that in 1 and 2 dimensions, rigidity of particle number follows from hyperuniformity and decay of the truncated pair correlation function (decay like $r^{-2}$ or faster in 1D and faster than $r^{-4}$ in 2D). This result covers all known examples of number rigidity in 1 and 2 dimensional particle systems. Apart from the previous examples, it also includes the 1D Dyson log gas at inverse temerature $\beta \le 2$ and Coulomb systems for small $\beta$ in dimension $d \ge 2$. It also includes, by the remark following \eqref{suphom1}, all processes in 1D that exhibit a bounded variance of particle number, and perturbed lattice systems in 1 and 2 dimensions.

In any determinantal process, all statistical information is, in principle, encoded in the pair $(K,\mu)$, where $K$ is the kernel and $\mu$ is the background measure. In view of this, it is a pertinent question as to whether we can read off any aspect of the rigidity behaviour of the process by testing some simple properties of the pair $(K,\mu)$. In this direction, it has been shown in \cite{GK} that, in any general determinantal process (not necessarily on a Euclidean space), there is number rigidity \textsl{only if} $K$ is the kernel of an integral operator that acts as a projection on $L^2(\mu)$. This is consistent with the conjecture that hyperuniformity is a necessary condition for rigidity.

%In the case of i.i.d. perturbations of a lattice, it is known that in 1 and 2 dimensions, 
\cite{PS} investigated the rigidity of the i.i.d. perturbation of $\Z^d$. For $d=1,2$, they showed that there is rigidity of numbers as soon as the random perturbation has a finite $d$-th moment. This is consistent with the results of \cite{GL}. For $\Z^d, d>2$, \cite{PS}  showed that for Gaussian perturbations  there is a phase transition in the rigidity behaviour in terms of the standard deviation $\sigma$ of the Gaussian. When $\sigma$ is below a critical $\sigma_c$, there is number rigidity, and when $\sigma>\sigma_c$, there is no rigidity.  This, in particular, negates any possibility for a sufficiency criterion for number rigidity (on the lines of \cite{GL}) in dimensions $d>2$, since for the Gaussian perturbation the truncated pair correlation decays exponentially for all $\sigma$ (as shown by the formulae in Section \ref{lattice}).

% To see the connection between super-homogeneity and rigidity, we note that the variance of a linear statistic $\Var \left(\sum_{x \in \text{point process}} \phi(x)\right)=V(\phi)$ can be written as \[V(\phi)=\int_{\R^d} \int_{\R^d} \phi(x) \phi(y) G(x,y)dxdy.\]

\section{Higher rigidity}
The plethora of highly interesting instances in nature of the phenomenon of \textsl{number rigidity} naturally raises the question as to whether there are other manifestations of such \textsl{rigidity phenomena}, particularly involving statistics other than a simple particle count. The first result in this direction was obtained in \cite{GP}, where it was shown that in the standard planar Gaussian zero process, for any bounded open set $\L$, the point configuration $X_{\L^\c}$ outside $\L$ determines precisely the number and the sum of the points inside $\L$ (equivalently, the mass and the centre of mass of the particles in $\L$). It was further established that, subject to the constraint on the number and the sum (imposed by the configuration outside), the particles inside $\L$ could be in any \textsl{generic} location inside $\L$ with positive probability density (with respect to the Lebesgue measure on the relevant conserved sub-manifold).

Subsequently, this result has been widely generalized in \cite{GK} to the case of $\a$-Gaussian zeros. In particular, it has been shown that for the zeros of the $\a$-Gaussian entire functions, for any bounded open set $\L$, the outside configuration $X_{\L^\c}$ almost surely determines the first $\l(\lfloor \frac{1}{\a} \rfloor +1\r)$ (holomorphic) moments of the points inside $\L$. Furthermore, subject to these constraints, the inside points could be in any \textsl{generic} configuration inside $\L$ with positive probability density (with respect to the appropriate Lebesgue measure).

\section{Proof techniques}
\label{prooftech}
The basic idea of \cite{GP}, \cite{G} and \cite{GL} to prove number rigidity of a point process $\Xi$ is to find a sequence of functions $\phi^{[\eps]}(x)$ such that, $\phi^{[\eps]}(x)=1$ for $x\in \Lambda$ and \[ \var\l(\sum_{x_i \in \Xi}\phi^{[\eps]}(x_i)\r) \le \epsilon,\] for any $\epsilon>0$. Then for small  $\epsilon\to0$ we have 
\[ \sum_{x_i \in \Xi}\phi^{[\eps]}(x_i)  =\sum\chi_{\Lambda}(x_i)+\sum\chi_{\Lambda^c}(x_i)\phi^{[\eps]}(x_i)\]
\[=N(\Lambda)+\sum\chi_{\Lambda^c}(x_i)\phi^{[\eps]}(x_i)\]
\[\sim \E \l[ \sum \phi^{[\eps]}(x_i)\r]=\int \rho(x)\phi^{[\eps]}(x)dx,\] 
where $\chi_{\Lambda}(x)$ is the characteristic function of the set $\Lambda$. This determines $N_{\Lambda}$ given $X_{\Lambda^c}$.

This is accomplished in the most basic cases by choosing a sequence $\phi_R(x)=\phi(x/R)$ with an appropriate $\phi(x)$. More sophisticated situations demand a Cesaro-type mean of a number of such functions in order to achieve  the low-variance criterion.

To give a concrete example, we consider the case of number rigidity for the zeroes of the standard planar GAF (i.e., Gaussian Analytic Function). In this case, it is known that, if $\ph$ is a $C_c^2$ function and $\ph_L(\cdot):= \ph(\cdot/L)$, then 
\begin{equation}
\label{gafasy}
\var \l(  \sum_{x_i \in \Xi}\ph_L(x_i)  \r)  \xrightarrow{ L \to \infty } C \| \Delta \ph \|_2^2/L^2.
\end{equation} 
Thus, to prove number rigidity for $\L$ the unit disk following the approach mentioned above, we choose $\Phi$ to be a $C_c^2$ function that is $\equiv 1$ on $\L$, and $L$ to be large enough (depending on $\eps$, such that $\var \l(  \sum_{x_i \in \Xi}\Phi_L(x_i)  \r) \le \eps$, which can accomplished due to \eqref{gafasy}).

On the other hand, for the Ginibre ensemble, it is known that for $\ph$ and $\ph_L$ defined as in \eqref{gafasy} we have
\begin{equation}
\label{ginasy}
\var \l(  \sum_{x_i \in \Xi}\ph_L(x_i)  \r)  \xrightarrow{ L \to \infty } C \| \nabla \ph \|_2^2.
\end{equation}
Due to this, a choice of $\Phi$ similar to the GAF case cannot be made directly. To overcome this difficulty, we consider a $C_c^2$ function $\phi$ that is $\equiv 1$ on $\L$, and look at the various scaling $\phi_{2^n}$ of $\phi$. For $L=2^N$, we then define 
\[\Phi^{[N]}:=\frac{1}{N} \l( \sum_{j=1}^N \phi_{2^j} \r).\]
This is an analogue of a Cesaro-type sum of the various scalings $\phi_{2^j}$ of $\phi$. It can be shown that the random sums $\l(  \sum_{x_i \in \Xi}\ph_{2^j}(x_i)  \r)_{j=1}^{\infty} $ exhibit a fast decay of correlations at widely different scales $2^j,2^k$. This can be used to show that $\var(\Phi^{[N]}) \to 0$ as $N \to \infty$, and the rest of the proof can then be completed as in the case of the GAF zeros. 

\section{Outlook}
\label{outlook}
In \cite{GL}, the authors provide sufficient criteria for number rigidity in dimensions 1 and 2, in terms of hyperuniformity and decay of correlations. It is an intriguing question to ask whether hyperuniformity, along with appropriate assumptions on the decay of correlations, are in fact necessary for rigidity phenomena. Such a conjecture is in a sense supported by the following big-picture heuristic. When $\var[N(\L)]$ grows like $|\L|$,  (to the leading order) it behaves like an additive functional on two adjacent domains. This appears to indicate that surface effects become inconsequential in the limit $|\L| \to \infty$, which does not seem to be consistent with number rigidity. It is a pertinent question to explore whether such criteria can be found in dimensions $d \ge 3$.
%Another direction in which little is known involves rigidity phenomena in dimensions $d \ge 3$. In particular, the ``spectral approach'', outlined in Section \ref{prooftech}, seems to be rather ineffective in higher dimensions, and suggests the necessity of novel methodologies. In \cite{GK}, the authors provide a hierarchy of rigidity behaviour in 2D for zeros of a 1-parameter family of Gaussian analytic functions. It is unknown whether one can obtain a  similar range of rigidity phenomena under the additional constraint that the point process be translation invariant. 

In \cite{G2}, the author makes a connection between rigidity phenomena and mutual regularity and singularity properties of Palm measures for very general point processes. E.g., for the zeros of the standard planar GAF, it is shown that the Palm measures at two points $z,w \in \C$, denoted resp. $\P_z,\P_w$, are mutually singular for Lebesgue a.e.-pair $(z,w)$. It is an interesting question to ask if this can be extended to cover all pairs $(z,w)$ with $z \ne w$, and if not, what is a description of the exceptional pairs? On a broader scale, it is pertinent to ask similar questions for mutual singularity of Palm measures in the generality considered in \cite{G2}.

\section{Acknowledgements}
The work of J.L.L. was supported in part by the NSF grant DMR1104501 and the AFOSR grant FA9550-16-1-0037. The work of S.G. was supported in part by the ARO grant W911NF-14-1-0094.

\begin{tabular}{l r}
\textsc{Subhroshekhar Ghosh}  \hspace{60 pt} & \textsc{Joel L. Lebowitz} \\
\textsc{Dept of ORFE}  & \textsc{Depts of Mathematics and Physics} \\
\textsc{Princeton University} & \textsc{Rutgers University} \\
% \textsc{Princeton, NJ 08544} & \textsc{Bangalore 560 012} \\
% \textsc{USA} & \textsc{India} \\
email: subhrowork@gmail.com & email: lebowitz@math.rutgers.edu \\ 
\end{tabular}

\end{document}